\documentclass[12pt]{article}
\usepackage{graphicx}
\usepackage[font=footnotesize,labelfont=bf,margin=4em]{caption}

\usepackage{amsfonts}
\usepackage{bbm}
\usepackage{mathrsfs}
 \usepackage{mathrsfs,amsmath,amssymb,amsthm}

\usepackage[dvips]{color}
 \usepackage{amssymb}
 \usepackage{amsbsy}
 \usepackage{color}
\allowdisplaybreaks
 \theoremstyle{plain}
\theoremstyle{remark}  \newtheorem{remark}{\noindent\mbox{Remark}}
 \theoremstyle{plain}
 \theoremstyle{plain}\newtheorem{lemma}{\noindent\mbox{Lemma}}
\theoremstyle{plain} \newtheorem{theorem}{\noindent\mbox{Theorem}}
 \theoremstyle{plain}\newtheorem{proposition}{\noindent\mbox{Proposition}}
 \theoremstyle{plain}\newtheorem{corollary}{\noindent\mbox{Corollary}}
\theoremstyle{definition} \newtheorem{definition}{\noindent\mbox{Definition}}

 \def\bq{\begin{equation}}
 \def\eq{\end{equation}}
 \def\eqn{\end{eqnarray}}
 \def\bqn{\begin{eqnarray}}
 \def\proof{\noindent{\it Proof.~~}}
 \def\qed{\hfill$\Box$\medskip}
 \def\rto{\rightarrow\infty}
 \def\z{\left}
 \def\y{\right}
 \def\no{\nonumber}

\begin{document}
 \title{\textbf{
 Maximum likelihood estimator and its consistency for an $(L,1)$ random walk in a parametric random environment}\footnote{Supported by
Nature Science Foundation of China(Grant No. 11501008;11801596).}}                  

\author{  Hua-Ming \uppercase{Wang}\footnote{Email:hming@ahnu.edu.cn; Department of Statistics, Anhui Normal University, Wuhu 241003, China }~$\ \&$ Meijuan \uppercase{Zhang}\footnote{Email: zhangmeijuan1227@163.com; School of Statistics and Mathematics, Central University of Finance and Economics, Beijing 100081, China }  }
\date{}
\maketitle

\vspace{-.8cm}

\begin{center}
\begin{minipage}[c]{12cm}
\begin{center}\textbf{Abstract}\quad \end{center}
Consider an $(L,1)$ random walk in an i.i.d. random environment, whose environment involves certain parameter.  We get the maximum likelihood estimator(MLE) of the environment parameter which can be written as functionals of a multitype branching process with immigration in a random environment(BPIRE). Because the offspring distributions of the involved multitype BPIRE are of the linear fractional type, the limit invariant distribution of the multitype BPIRE can be computed explicitly. As a result, we get the consistency of the MLE. Our result is a generalization of Comets et al. [Stochastic Process. Appl. 2014, 124, 268-288].
\vspace{0.2cm}

\vspace{0.2cm}

\end{minipage}
\end{center}

\section{Introduction}

Random walks in random environments(RWRE) exhibit many  surprising phenomena, thus attracting much attention recent years. Their limit behaviors, especially for one-dimensional case, have been extensively studied. We refer the readers to \cite{ze04} for a survey. RWREs involve two kinds of randomness: the environments are chosen randomly according to some distribution; the particle evolves randomly in a given environment. So, from the statistical point of view, it is interesting to infer the distribution of the environment up to a single observation of a path of the random
walk until it reaches a distant site. Adelman and Enriquez \cite{ae} dealt with very general RWRE and presented a procedure to infer
the environment distribution through a system of moment equations. For the specific nearest neighbor ballistic RWRE,  Comets et al. \cite{cfl} provided  a maximum
likelihood estimator(MLE) for the parameter of the environment distribution and studied its consistency, whereas Falconnet et al. \cite{FL14} studied the asymptotic normality and efficiency of the MLE.
It turns out that the MLE  of Comet et al. exhibits a
much smaller variance than the one of Adelman and Enriquez.

In this paper, we study a non-nearest neighbor RWRE, say $(L,1)$ RWRE, whose biggest left-oriented jumps are of size $L$ and right-oriented jumps are always with size $1.$ Motivated by Comets et al. \cite{cfl}, our main goal is to provide a procedure to give the MLE of the environment parameter and show its consistency. We adopt the approach used in \cite{cfl}. The MLE can be written as functionals of a multitype branching process with immigration in a random environment(BPIRE), which we get by decomposing the path of the RWRE. It turns out that the offspring distributions of the  involved multitype BPIRE are of the  linear fractional type. So we can compute explicitly the distribution of that BPIRE and its limit invariant distribution as well. Consequently, we can get the consistency of the MLE.

{\it Notes:} For general $L\ge2,$ the construction of the  MLE for $(L,1)$ RWRE is basically the same as (2,1) RWRE, except that for $(L,1)$ RWRE, a $(1+2+...+L)$-type BPIRE is involved while for $(2,1)$ RWRE, a $(1+2)$-type BPIRE is needed. However, for general $(L,1)$ RWRE, the notations are very heavy. So in the remainder of the paper, we fix $L=2$ to consider $(2,1)$ RWRE.

The paper is organized as follows. We devote Section \ref{mainre} to introducing the model and giving a procedure to construct the MLE. In Section \ref{br} we construct a 3-type BPIRE by decomposing the path of $(2,1)$ RWRE and compute  explicitly its limit distribution. Finally, in Section \ref{consi} we show that the MLE we construct is consistent.

\section{Model and main results}\label{mainre}
\subsection{  (2,1) RWRE and its Preliminaries}
For $x\in \mathbb Z,$ let $\omega_x=(\omega_x(-2),\omega_x(-1),\omega_x(1))$ be a probability measure on $\{x-2,x-1,x+1\}$.
Let $\Omega$ be the collection of all $\omega=(\omega_x,~x\in\mathbb Z).$ Equip $\Omega$ with the weak topology induced by convergence of probability measures and let $\mathcal F$ be the Borel $\sigma$-algebra. Let $\nu_\theta$ be the law of $\omega_0$ with  $\theta\in\Theta$ certain unknown parameter. We always assume that $\Theta\subset\mathbb{R}^{d}$ for some  $d\ge 1$ is a compact set.  Then $\mathbb{P}^{\theta}:=\nu_{\theta}^{\otimes\mathbb{Z}}$ is a probability measure on $(\Omega, \mathcal F)$ which makes $\omega=(\omega_x,~x\in\mathbb Z)$ an i.i.d. sequence. For a realization of $\omega,$ we consider a Markov chain $\{X_t\}_{t\ge 0}$ on $\mathbb{Z}$ starting from $0$, with transition probabilities
\begin{equation*}
P_{\omega}(X_{t+1}=x+l\big |X_{t}=x)=\omega_{x}(l),\text{ for } l=-2,-1,1
\end{equation*}
so that $P_\omega$ is the quenched law of the Markov chain under the
environment $\omega.$ Define a new probability measure $\mathbf{P}^{\theta}$ by $$\mathbf{P}^{\theta}(\cdot)=\int P_\omega(\cdot)\mathbb{P}^{\theta}(d\omega),$$ which is usually called the annealed probability of $\{X_t\}.$ We denote by $\mathbb{E}^{\theta}$, $E_\omega$ and $\mathbf{E}^{\theta}$ the expectation operators for $\mathbb{P}^{\theta}$, $P_\omega$ and $\mathbf{P}^{\theta}$, respectively.
To state the recurrence criterion, we need the following condition.

\bigskip
\noindent{\bf (C1)} {\it Suppose that
 $\inf_{\theta\in\Theta}\mathbb{E}^{\theta}(\log\omega_0(l))>-\infty,$ for $ l=-2,-1,1.$}

\bigskip
For $k\in\mathbb{Z},$ let
 $a_{k}=\frac{\omega_{k}(-1)}{\omega_{k}(1)},$  $b_{k}=\frac{\omega_{k}(-2)}{\omega_{k}(1)}$ and set
\begin{equation}\label{ad}
A_{k}=\left(
\begin{array}{ccc}
a_k & b_k & 0 \\
a_k & b_k  & 1 \\
a_k & b_k &  0 \\
\end{array}
\right),\ B_k
=\left(
\begin{array}{cc}
a_k+b_k&b_k \\
1 & 0 \\
\end{array}
\right).
\end{equation}
Clearly,  $\{A_k\}$ and $\{B_k\}$ are  two sequences of i.i.d. random matrices under $\mathbb{P}^{\theta}.$  Then (C1) ensures an application of Oseledec's multiplicative ergodic theorem (see \cite{osel}) to obtain the Lyapunov exponents of $\{A_k\}$ and $\{B_k\}$  under $\mathbb P^{\theta}.$ Let $\gamma_A$ and $\gamma_B$ be the top Lyapunov exponents of $\{A_k\}$ and $\{B_k\},$ respectively. Then by the positivity of entries in $A_k$ and $B_k,$ we have (see \cite{king}, Theorem 5) that $\mathbb P^{\theta}$-a.s., for all $i,j\in \{1,2,3\},$
\begin{equation}\label{gab}
\gamma_A=\lim_{n\rightarrow\infty}\frac{1}{n}\log \mathbf{e}_iA_0
A_1\cdots A_{n-1}\mathbf{e}_j^{t},\ \gamma_B=\lim_{n\rightarrow\infty}\frac{1}{n}\log \mathbf{e}_iB_0
B_1\cdots B_{n-1}\mathbf{e}_j^{t}.
\end{equation}
  Here and throughout, $\mathbf e_i$ is a vector with the $i$-th component $1$ and all other components $0,$ and $\mathbf v^{t}$ denotes the transpose of a vector $\mathbf v.$ Unless otherwise stated, the dimension of a vector depends on the matrix multiplication.
By induction, it can be easily verified that
$\mathbf{e}_{1} A_{0}A_{1}\cdots A_{k-1}\mathbf{e}_{2}^{t}
=\mathbf{{e}}_{1}B_{0}B_{1}\cdots B_{k-1}\mathbf{{e}}_{2}^{t}.
$
Thus we can infer from (\ref{gab}) that $$\gamma_A=\gamma_B.$$
The following recurrence criterion can be found in  Letchikov \cite{letca}.

{\noindent \bf Recurrence criterion:} {\it Under (C1),  $\{X_t\}$ is transient to $+\infty,$ recurrent or transient to $-\infty$ according as $\gamma_A<0,$ $\gamma_A=0$ or $\gamma_A>0,$ respectively.}

\begin{remark}
  Here,  $\gamma_A,$ the top Lyapunov exponent of $\{A_k\},$ is used to give the recurrence criterion whereas in \cite{letca},  $\gamma_B$ for $\{B_k\}$ is used. This causes no problem since we have shown that $\gamma_A=\gamma_B.$ We use $\gamma_A$ here because we will see below that $\{A_k\}$ serves as the quenched offspring mean matrices for a $3$-type BPIRE constructed from the path of $\{X_t\}.$
\end{remark}
In this paper, we consider the case that $\{X_t\}$ is transient to $+\infty.$  That is, we assume the following condition holds.
\bigskip

\noindent{\bf (C2)} {\it Suppose that $\gamma_A<0.$}

\bigskip
Next we give a condition to ensure that $\{X_t\}$ has a positive speed or in other words, $\{X_t\}$ is ballistic. Let $$\pi(\omega)=1+\sum_{k=1}^\infty \mathbf e_1A_k\cdots A_1(2,1,2)^t.$$
{\noindent \textbf{(C3)} {\it Suppose that  for any $\theta\in\Theta,$ we have $\mathbb E^\theta (\pi(\omega))<\infty.$}

\bigskip
\noindent {\bf Law of large numbers:} {\it Suppose that (C1)-(C3) hold. Then for any  $\theta\in \Theta,$ we have $\mathbf P^{\theta}$-a.s., $\lim_{t\rto}\frac{X_t}{t}=\frac{1}{\mathbb E^\theta(\pi(\omega))}.$}

Let $T_1:=\inf\{t>0: X_t=1\}.$ Then under (C1)-(C3), from Corollary \ref{mm} below and the stationarity of the environment, we know that $\mathbf E^\theta T_1=\mathbb E^\theta(\pi(\omega))<\infty.$ Consequently, we can mimic the proof of Theorem 2.1.9 in \cite{ze04} to prove the above law of large numbers by a hitting time decomposition approach.

\subsection{Construction of an $M$-estimator}

For $n\ge 1,$ define $$T_n:=\inf\{t>0:X_t=n\},$$ the time the walk hits  $n$ for the first time. Under (C2),  $\lim_{t\rto} X_t=+\infty,$ $\mathbf{P}^{\theta}$-a.s. and hence for each $n\ge1,$ $T_{n}<\infty,$ $\mathbf{P}^{\theta}$-a.s. For  $n\geq1,$ an  observation of the path of $\{X_t\}$ until it reaches $n$ is denoted by $X_{[0,T_{n}]}=\{X_{t}:~t=0,1,\cdots,T_{n}\}.$

Let $\mathbf{x}_{[0,t]}:=(x_{0},\cdots,x_{t})$ be a path of length $t,$  $v_t$ be the set of integers visited by the path $\mathbf{x}_{[0,t]}$  and $V_{T_{n}}$ be the set of integers visited by the path $X_{[0,T_{n}]}$. For $x\in\mathbb{Z},$  define
\begin{align}
&L_{1}(x,\mathbf{x}_{[0,t]})=\#\{0\le s<t: x_{s}=x,x_{s+1}=x-1\},\no\\
&L_{2}(x,\mathbf{x}_{[0,t]})=\#\{0\le s<t: x_{s}=x,x_{s+1}=x-2\},\no\\
&R(x,\mathbf{x}_{[0,t]})=\#\{0\le s<t: x_{s}=x,x_{s+1}=x+1\}.\no
\end{align}
Here and throughout, we denote by $\#\{\ \}$ the number of elements in a set $\{\ \}.$ Clearly, we have
\begin{align}
P_{\omega}(X_{[0,T_{n}]}=\mathbf{x}_{[0,t_{n}]})
&=\prod_{x\in v_{t_n}}\omega_{x}(1)^{R(x,\mathbf{x}_{[0,t_{n}]})}
\omega_{x}(-1)^{L_{1}(x,\mathbf{x}_{[0,t_{n}]})}
\omega_{x}(-2)^{L_{2}(x,\mathbf{x}_{[0,t_{n}]})},\no\\
\mathbf{P}^{\theta}(X_{[0,T_{n}]}=\mathbf{x}_{[0,t_{n}]})
&=\prod_{x\in v_{t_n}}~~\iint\limits_{0\leq a_{1}+a_{2}\leq 1}(1-a_{1}-a_{2})
^{R(x,\mathbf{x}_{[0,t_{n}]})}\no\\
&\quad\quad\quad\quad\quad\quad\quad\quad\times
a_{1}^{L_{1}(x,\mathbf{x}_{[0,t_{n}]})}
a_{2}^{L_{2}(x,\mathbf{x}_{[0,t_{n}]})}d\nu_{\theta}(a_{1},a_{2}),\no
\end{align}
where $\mathbf{x}_{[0,t_{n}]}$ is the path up to time  $ t_{n},$  the first hitting time of site $n.$

One can define $L_{1}(x,X_{[0,T_{n}]}),$ $L_{2}(x,X_{[0,T_{n}]})$ and $R(x,X_{[0,T_{n}]})$ by a similar way.   Write simply
$
L_{x,1}^{n}:=L_{1}(x,X_{[0,T_{n}]}),$ $L_{x,2}^{n}:=L_{2}(x,X_{[0,T_{n}]}),$ $R_{x}^{n}:=R(x,X_{[0,T_{n}]})$ and set $$\mathbf{L}_{x}^{n}=(L_{x,1}^{n},L_{x,2}^{n}).$$
Obviously, $\mathbf{L}_{x}^{n}$ counts the left-oriented jumps at $x$ within the path $X_{[0,T_{n}]}.$

Next we define an annealed log-likelihood function based on  $X_{[0,T_{n}]},$ a single observation of the path of  $\{X_t\}$ until it reaches a distant site $n.$ For $n\ge 1,$ define
\begin{equation}\label{likeli2}\begin{split}
&\tilde{l}_{n}(\theta):=\sum_{x=0}^{n-1}\log\iint\limits_{0\leq a_{1}+a_{2}\leq 1}a_{1}^{L_{x,1}^{n}} a_{2}^{L_{x,2}^{n}}(1-a_{1}-a_{2})
^{R_{x}^{n}}d\nu_{\theta}(a_{1},a_{2})
\\
&\quad\quad\quad\quad+\displaystyle\sum_{x<0,~x\in V_{T_{n}}}\log\iint\limits_{0\leq a_{1}+a_{2}\leq1}a_{1}^{L_{x,1}^{n}}a_{2}^{L_{x,2}^{n}}(1-a_{1}-a_{2})
^{R_{x}^{n}}d\nu_{\theta}(a_{1},a_{2}).\end{split}
\end{equation}
Under (C2), $\lim_{t\rto}X_t=+\infty,$ $\mathbf P^{\theta}$-a.s. Thus $\{X_t\}$ can only visit at most finite negative sites. Hence the second term on the righthand of (\ref{likeli2})  will not affect significantly the  behavior of the normalized log-likelihood function, say $\tilde{l}_{n}(\theta)/n.$ So we keep only the first term which plays a dominant role.

If $\mathbf P^{\theta}$-a.s. $\lim_{t\rto}X_t=+\infty,$ then for $x\ge0,$ we have
\begin{eqnarray}
R_{x}^{n}=L_{x+1,1}^{n}+L_{x+1,2}^{n}+L_{x+2,2}^{n}+1.\no
\end{eqnarray}
Therefore, if we denote the first term on the righthand of (\ref{likeli2}) by $l_n(\theta),$ then
\begin{equation}\label{likeli}
l_{n}(\theta)
=\sum_{x=0}^{n-1}\log\iint\limits_{0\leq a_{1}+a_{2}\leq 1}a_{1}^{L_{x,1}^{n}} a_{2}^{L_{x,2}^{n}}
(1-a_{1}-a_{2})^{L_{x+1,1}^{n}+L_{x+1,2}^{n}+L_{x+2,2}^{n}+1}
d\nu_{\theta}(a_{1},a_{2}).\no
\end{equation}
$l_n(\theta)$ will serve as the criterion function for deriving an $M$-estimator for  $\theta.$ Set $\mathbb N=\{0,1,2,...\}$ and
define a function  $\phi_{\theta}:\mathbb{N}^{2}\times \mathbb{N}^{2}\times \mathbb{N}^{2} \mapsto \mathbb{R}$ by setting
\begin{equation}\label{dph}
\phi_{\theta}(\mathbf{y}_{1},\mathbf{y}_{2},\mathbf{y}_{3})
=\log\iint\limits_{0\leq a_{1}+a_{2}\leq 1}
a_{1}^{y_{1,1}} a_{2}^{y_{1,2}}
(1-a_{1}-a_{2})^{y_{2,1}+y_{2,2}
+y_{3,2}+1}d\nu_{\theta}(a_{1},a_{2}),
\end{equation}
where $\mathbf{y}_{i}=(y_{i,1},y_{i,2})\in \mathbb{N}^{2},~i=1,2,3$.
Then
\begin{equation}\label{lnt}
l_{n}(\theta)=\sum_{x=0}^{n-1}\phi_{\theta}
(\mathbf{L}_{x}^{n},\mathbf{L}_{x+1}^{n},\mathbf{L}_{x+2}^{n}).
\end{equation}
We need the following conditions in addition.

 \noindent{\bf (C4)} {\it The map $\theta\mapsto\nu_{\theta}$ is continuous on $\Theta$ with respect to the weak topology.}

   \noindent{\bf (C5)}  {\it $\forall (\theta,\theta')\in\Theta^{2},$ ~$\nu_{\theta}\neq\nu_{\theta'}$
~$\Leftrightarrow$~$\theta\neq\theta'.$}

Note that under (C4), since $\Theta$ is compact, it follows from (\ref{dph}) that for any $(\mathbf{y}_{1},\mathbf{y}_{2},\mathbf{y}_{3}),$ $\phi_{\theta}(\mathbf{y}_{1},\mathbf{y}_{2},\mathbf{y}_{3})$ is continuous in $\theta.$ Thus $l_{n}(\theta)$ is also a continuous function of $\theta$ and hence   it can achieve its maximum over the compact set $\Theta.$ Therefore, we can define an estimator $\hat{\theta}_{n}$ of $\theta$ as follows.
\begin{definition}\label{def}
An estimator $\hat{\theta}_{n}$ of $\theta$ is defined as a measurable choice
\begin{equation}\label{arg}
\hat{\theta}_{n}\in \displaystyle\mbox{Argmax}_{\theta\in\Theta}l_{n}(\theta).
\end{equation}
\end{definition}
Since $\hat{\theta}_{n}$ is a maximiser of the criterion function $l_{n}(\theta),$ it is an  $M$-estimator. Obviously, it is not necessarily unique. The criterion function $l_{n}(\theta)$ is not exactly log-likelihood. But the contribution of the negative sites which we dropped is not so important. So, we still call $\hat{\theta}_{n}$ a MLE.

Assume that the process $\{X_{t}\}$ is generated under the true parameter value $\theta^{\ast},$  an interior point of the parameter space $\Theta,$ which we want to estimate. For simplicity, we use $\mathbf{P}^{\ast}$ and $\mathbf{E}^{\ast}$ rather than $\mathbf{P}^{\theta^{\ast}}$ and $\mathbf{E}^{\theta^{\ast}},$ respectively.

The following theorem, whose proof will be postponed to Section \ref{consi}, gives the consistency of the estimator $\hat\theta_n.$
\begin{theorem}[Consistency] \label{consist} Suppose that (C1)-(C5) hold. Then
for any choice of $\hat{\theta}_{n}$ satisfying (\ref{arg}),  we have
$\lim_{n\rightarrow+\infty}\hat{\theta}_{n}=\theta^{\ast},$
in $\mathbf{P}^{\ast}$-probability.
\end{theorem}
\section{A $3$-type BPIRE in the path of (2,1) RWRE}\label{br}
Throughout  this section, we always assume that (C1) and (C2) hold,  that is, $\lim_{t\rto}X_t=+\infty,$ $\mathbf P^\theta$-a.s.
\subsection{Construction of BPIRE from the path of  RWRE}
In this section, we construct a $3$-type BPIRE from the path of $\{X_t\}.$ Hong and Wang \cite{[HW09]} revealed a $2$-type BPIRE from (2,1) RWRE. But that result is not sufficient  to show the consistency for MLE. Here, we give a new construction, which is slightly different from  \cite{[HW09]}, but with the same idea.

Fix an integer $n>0.$ Let
\begin{eqnarray}
&& U_{i,1}^n=\#\{0<k<T_n:X_{k-1}=i,X_k=i-1\},\no\\
&& U_{i,2}^n=\#\{0<k<T_n:X_{k-1}=i,X_k=i-2\},\no\\
&& U_{i,3}^n=\#\{0<k<T_n:X_{k-1}=i+1,X_k=i-1\}\no
\end{eqnarray}
and set
$$\mathbf U^n_i:=(U_{i,1}^n,U_{i,2}^n,U^n_{i,3}).$$  From the definition, we see that $U_{i,1}^n,U_{i,2}^n$ and $U^n_{i,3}$
count the number of steps by the walk  from $i$ to $i-1,$ $i$ to $i-2$ and $i+1$ to $i-1$ respectively in the time interval $[0,T_n].$
\begin{figure}[h]
  \begin{center}
  \includegraphics[scale=.8]{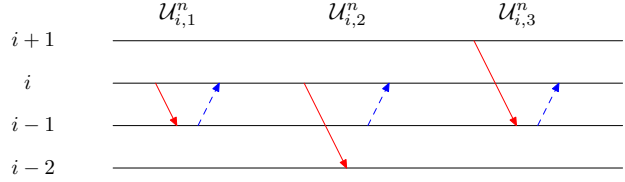}\\
  \caption{  The figure illustrates three types of  jumps at $i.$ Since $\lim_{t\rto}X_t=+\infty,$ $\mathbf P^\theta$-a.s., there must be an upward jump from $i$ to $i+1$ which matches the corresponding downward jump. So in the figure, an upward (dashed) jump from $i$ to $i+1$  is attached to each downward jump. Thus, we treat $\mathcal U_{i,l}$ an excursion rather than a jump, and such excursions constitute the individuals of a $3$-type BPIRE.}\label{ex}
    \end{center}
    \vspace{-1cm}

  \end{figure}

\begin{theorem} Assume (C1) and (C2). Then  $\mathbf U^n_n\equiv0, \mathbf U^n_{n-1}, \mathbf U^n_{n-2},...,  \mathbf U^n_{1}, \mathbf U^n_{0}$ form the first $n$ generations of a $3$-type BPIRE, which evolves as follows: given $\omega,$ at time $n-(i+1),$ $0\le i\le n-1$(in view of the walk, time $n-(i+1)$ of BPIRE corresponds to the site $i+1$), a type-$1$ individual immigrates into the system, which, together with the existing individuals that constitute $\mathbf U^n_{i+1},$ will give birth to a number of descendants with the law
  \begin{align}
  &P_\omega(\mathbf U^n_{i}=(a,b,0)|\mathbf U^n_{i+1}=\mathbf e_1)
=\frac{(a+b)!}{a!b!}(\omega_{i}(-1))^{a}(\omega_{i}(-2))^{b}\omega_{i}(1),\label{disa}\\
&P_\omega(\mathbf U^n_{i}=(a,b,1)|\mathbf  U^n_{i+1}=\mathbf e_2)
=\frac{(a+b)!}{a!b!}(\omega_{i}(-1))^{a}(\omega_{i}(-2))^{b}\omega_{i}(1),\label{disb}\\
&P_\omega(\mathbf U^n_{i}=(a,b,0)|\mathbf U^n_{i+1}=\mathbf e_3)
=\frac{(a+b)!}{a!b!}(\omega_{i}(-1))^{a}(\omega_{i}(-2))^{b}\omega_{i}(1).\label{disc}
\end{align}
Moreover, we have
\begin{equation}\label{tu}T_n=n+\sum_{k=-\infty}^{n-1}\mathbf U^n_k(2,1,2)^t.\end{equation}
\end{theorem}
\proof We give here only the idea of the proof. For more details, we refer the readers to \cite{[HW09]}. Given $\omega,$  various independence for the branching process follows from  the independence among different excursions, which follows from the strong Markov properties. As for the immigration law, note that for $1\le k\le n,$ before $\{X_t\}$ hits $k$ for the first time, there might be a number of excursions $\mathcal U_{k-1,1}$ and $\mathcal U_{k-1,2}.$ But before time $T_k,$ there is no  step that reaches $k-1$ from some site above $k-1.$ So we may treat such excursions $\mathcal U_{k-1,1}$ and $\mathcal U_{k-1,2}$ as individuals born to a type-$1$ individual, which immigrates into the system at site $k.$
For the offspring distributions, we only explain that of a type-$2$ individual, which we illustrate in Figure \ref{ot}.
\begin{figure}[h]
  \centering
  \includegraphics[scale=0.9]{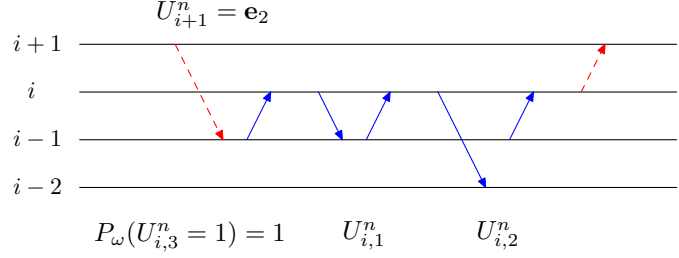}\\
  \caption{ The offspring distribution of a type-$2$ individual at $i+1.$}\label{ot}

  \vspace{-.4cm}

\end{figure}

We see from Figure \ref{ot} that given $\omega,$ a type-$2$ individual at $i+1$(the excursion starts from the dashed downward jump and ends up with the dashed upward jump)  gives birth to type-$3$ individual at $i$ with probability $1$ firstly and then  to a number of type-$1$ and type-$2$ individuals. If it produces $a$ type-$1$ and $b$ type-$2$ individuals, that is, $U_{i,1}=a,$ $U_{i,2}=b,$ then during the excursion at $i+1$ considered here, the walk will firstly jump  from $i$ to $i-1$ with $a$ times, jump from $i$ to $i-2$ with $b$ times and finally jump from $i$ to $i+1.$ All of the possible combinations of those $a+b$ jumps equal $\frac{(a+b)!}{a!b!}.$ We thus conclude that
$$P_\omega(\mathbf U^n_{i}=(a,b,1)|\mathbf U^n_{i+1}=\mathbf e_2)
=\frac{(a+b)!}{a!b!}(\omega_{i}(-1))^{a}(\omega_{i}(-2))^{b}\omega_{i}(1).$$
The offspring distribution for a type-$1$(or type-$3$) individual is much easier and can be discussed similarly. Also, by some delicate computation, we see that the righthand of (\ref{tu}) counts all steps before the walk hits $n.$ So (\ref{tu}) is true. \qed

From (\ref{disa})-(\ref{disc}), we get the mean offspring matrix as follows.
\begin{corollary}\label{mm} Given $\omega,$ let $A_i$ be a $3\times 3$ matrix whose $l$-th row is $E_\omega(U_{i}^n|U_{i+1}^n=\mathbf e_l).$ Then $A_i$ coincides with the one defined in (\ref{ad}). Moreover
\begin{align*}
  &E_\omega T_n=n+\sum_{i=0}^{n-1}\sum_{k=-\infty}^{i}\mathbf e_1 A_{i}\cdots A_{k+1}A_{k}(2,1,2)^t.
\end{align*}
 \end{corollary}

%
%
%
%

Next, we define a new  $3$-type BPIRE $\{\mathbf Z_n\}_{n\ge0}.$
  Firstly we specify its offspring distributions. For $1\le i\le 3,$ $k\ge 1,$ let  $$\mathbf{\xi}(i,k)=(\xi_{1}(i,k),\xi_{2}(i,k),\xi_{3}(i,k)),$$ which will serve as offsprings in $k$-th generation born to a type-$i$ individual in $(k-1)$-th generation, be a random vector satisfying
\begin{eqnarray}
&P_\omega(\mathbf{\xi}(1,k)=(a,b,0))
=\frac{(a+b)!}{a!b!}(\omega_{k}(-1))^{a}(\omega_{k}(-2))^{b}
\omega_{k}(1),\no\\
&P_\omega(\mathbf{\xi}(2,k)=(a,b,1))
=\frac{(a+b)!}{a!b!}(\omega_{k}(-1))^{a}(\omega_{k}(-2))^{b}
\omega_{k}(1),\no\\
&P_\omega(\mathbf{\xi}(3,k)=(a,b,0))
=\frac{(a+b)!}{a!b!}(\omega_{k}(-1))^{a}(\omega_{k}(-2))^{b}
\omega_{k}(1).\no
\end{eqnarray}
Moreover, suppose that $\mathbf{\xi}(i,k),$ $i=1,2,3,$ $k\ge 1$ are mutually independent.

Set $\mathbf Z_0=\mathbf 0$ and for $n\ge 1,$ define recursively
\begin{align}\label{zr}
 \mathbf Z_n=\sum_{j=1}^{1+Z_{n-1,1}} \xi^j(1,n) +\sum_{j=1}^{Z_{n-1,2}} \xi^j(2,n) +\sum_{j=1}^{Z_{n-1,3}} \xi^j(3,n),\no
\end{align}
where for fixed $i=1,2,3$ and $n\ge 0,$ $\xi^j(i,n),j=1,2,...$ are independent copies of $\xi(i,n).$
Clearly, $\{\mathbf Z_n\}_{n\ge0}$ is a 3-type branching process with exactly one type-$1$ immigrant in each generation in a random environment.
Since $(\omega_{n-1},\cdots,\omega_{0})$ and $(\omega_{1},\cdots,\omega_{n})$ share the same law
 under $\mathbb{P}^\theta,$  we have
 \begin{equation}\label{uz}\begin{split}
  &\text{under }\mathbf P^\theta, \mathbf U^n_{n}\equiv 0,\mathbf U^{n}_{n-1},\cdots,
\mathbf U_{0}^n\text{ has the same distribution }\\
&\text{as }\mathbf{Z}_{0},\mathbf{Z}_{1},\cdots,\mathbf{Z}_{n},\text{ the first }n\text{ generations of }\{\mathbf Z_n\}.
\end{split}
\end{equation}
In view of (\ref{uz}),  we need only work with $\{\mathbf Z_n\}.$

\subsection{Criterion function $l_n(\theta)$ as functional of $\{\mathbf Z_n\}$}

It is easily seen that
\begin{equation}
U_{x,1}^n=L_{x,1}^n,\quad U_{x,2}^n=L_{x,2}^n,\quad
U^n_{x,3}=U^n_{x+1,2}=L_{x+1,2}^n.\no
\end{equation}
Thus, we have
\begin{align}
R_{x}^{n}&=L_{x+1,1}^{n}+L_{x+1,2}^{n}+L_{x+2,2}^{n}+1\no\\
&=U_{x+1,1}^{n}+U_{x+1,2}^{n}+U_{x+1,3}^{n}+1=|\mathbf{U}^n_{x+1}|+1.\no
\end{align}
Consequently,
\begin{equation}
l_{n}(\theta)
=\sum_{x=0}^{n-1}\log\iint\limits_{0\leq a_{1}+a_{2}\leq 1}a_{1}^{U_{x,1}^{n}} a_{2}^{U_{x,2}^{n}}
(1-a_{1}-a_{2})^{|\mathbf{U}^n_{x+1}|+1}
d\nu_{\theta}(a_{1},a_{2}).\no
\end{equation}
Let $\tilde{\phi}_{\theta}$ be a function from $\mathbb{N}^{3}\times \mathbb{N}^{3}$ to  $\mathbb{R}$ defined by
\begin{equation}
\tilde{\phi}_{\theta}(\mathbf{z}_{1},\mathbf{z}_{2})
=\log\iint\limits_{0\leq a_{1}+a_{2}\leq 1}
a_{1}^{z_{2,1}}a_{2}^{z_{2,2}}(1-a_{1}-a_{2})^{z_{1,1}+z_{1,2}+z_{1,3}+1}
d\nu_{\theta}(a_{1},a_{2}),\no
\end{equation}
where $\mathbf{z}_{i}=(z_{i,1},z_{i,1},z_{i,3})\in\mathbb{N}^{3},~i=1,2.$
Then from (\ref{lnt}), we have
\begin{equation}
l_{n}(\theta)=\sum_{x=0}^{n-1}\phi_{\theta}
(\mathbf{L}_{x}^{n},\mathbf{L}_{x+1}^{n},\mathbf{L}_{x+2}^{n})
=\sum_{x=0}^{n-1}\tilde{\phi}_{\theta}
(\mathbf U_{x+1}^{n},\mathbf U^n_{x}).\no
\end{equation}
Therefore, it follows from (\ref{uz}) that
\begin{equation}\label{lz}
l_{n}(\theta)=
\sum_{k=0}^{n-1}\tilde{\phi}_{\theta}
(\mathbf{Z}_{k},\mathbf{Z}_{k+1})
\end{equation}
in $\mathbf{P}^{\theta}$-distribution.
Thus, to  characterize  $l_{n}(\theta),$ it is enough to study the properties of $\{\mathbf{Z}_{n}\}.$

\subsection{Quenched probability generating function of $\{\mathbf Z_n\}$}\label{mbp}

Given $\omega,$ since the offspring distributions of $\{\mathbf Z_n\}$ are of the linear fractional type, we can compute explicitly its probability generating function(p.g.f.).
%

As a convention, set $\mathbf{s}=(s_{1},s_{2},s_{3})^{t}\in[0,1]^{3}$ and $\mathbf{1}=(1,1,1)^{t}.$
For $k\ge0,$  denote by \begin{equation}
F_{k}^{1}(\mathbf{s})=E_{\omega}\big(\mathbf{s}^{\mathbf{Z}_{k}}\big)
=E_{\omega}\big(s_{1}^{Z_{k,1}}s_{2}^{Z_{k,2}}s_{3}^{Z_{k,3}}
\big)\no
\end{equation}
the p.g.f. of $\mathbf Z_k.$  To give a formula for $F_{k}^{1}(\mathbf{s}),$ we introduce
\begin{align}
&S_{n}^{1}=\mathbf{e}_{1}\sum_{j=1}^{n}\prod_{i=j}^{n}A_{i}\mathbf{e}_{1}^{t},\
S_{n}^{2}=\mathbf{e}_{1}\sum_{j=1}^{n}\prod_{i=j}^{n}A_{i}\mathbf{e}_{2}^{t},\
S_{n}^{3}=\mathbf{e}_{1}\sum_{j=1}^{n}\prod_{i=j}^{n}A_{i}\mathbf{e}_{3}^{t},\no\\
&\widetilde{S}_{n}^{1}
=\mathbf{e}_{1}\sum_{j=1}^{n}\prod_{i=1}^{j}A_{i}\mathbf{e}_{1}^{t},\
\widetilde{S}_{n}^{2}
=\mathbf{e}_{1}\sum_{j=1}^{n}\prod_{i=1}^{j}A_{i}\mathbf{e}_{2}^{t},\
 \widetilde{S}_{n}^{3}
=\mathbf{e}_{1}\sum_{j=1}^{n}\prod_{i=1}^{j}A_{i}\mathbf{e}_{3}^{t},\no
\end{align}
and let $$S_{n}=1+S_{n}^{1}+S_{n}^{2}+S_{n}^{3},\ \widetilde{S}_{n}=1+\widetilde{S}_{n}^{1}+\widetilde{S}_{n}^{2}
+\widetilde{S}_{n}^{3}.$$
The following proposition is the main result of this section.
\begin{proposition}\label{Flem4}
 We have
\begin{equation}\label{ff}
F_{n}^{1}(\mathbf{s})
=\frac{\frac{1}{S_{n}}}{1-\frac{S_{n}^{1}}{S_{n}}s_{1}
-\frac{S_{n}^{2}}{S_{n}}s_{2}-\frac{S_{n}^{3}}{S_{n}}s_{3}},
\end{equation}
which equals
\begin{equation}\label{tff}\no
\widetilde{F}_{n}^{1}(\mathbf{s})
:=\frac{\frac{1}{\widetilde{S}_{n}}}{1-\frac{\widetilde{S}_{n}^{1}}
{\widetilde{S}_{n}}s_{1}
-\frac{\widetilde{S}_{n}^{2}}{\widetilde{S}_{n}}s_{2}
-\frac{\widetilde{S}_{n}^{3}}{\widetilde{S}_{n}}s_{3}}
\end{equation}
in $\mathbb P^\theta$-distribution.
\end{proposition}
To prove Proposition \ref{Flem4}, we need several lemmas for preparation.
Recall that for $i=1,2,3$ and $k\ge 0,$  $\mathbf{\xi}(i,k)=(\xi_{1}(i,k),\xi_{2}(i,k),\xi_{3}(i,k))$ are offsprings in $k$-th generation born to a type-$i$ individual in $(k-1)$-th generation. Given $\omega,$ for $k\ge0$ and $i\in\{1,2,3\}$ let
 \begin{equation}
f_{k}^{i}(\mathbf{s})
=E_{\omega}(\mathbf{s}^{\mathbf{\xi}(i,k)})
=E_{\omega}(s_{1}^{\xi_{1}(i,k)}s_{2}^{\xi_{2}(i,k)}s_{3}^{\xi_{3}(i,k)})
\no
\end{equation}
be the p.g.f. of $\mathbf{\xi}(i,k)$ and set $$\mathbf{f}_{k}(\mathbf{s})=(f_{k}^{1}(\mathbf{s}), f_{k}^{2}(\mathbf{s}),f_{k}^{3}(\mathbf{s}))^{t}.$$

\begin{lemma}\label{Flem1}
 For $k\geq 1,$ we have
$F_{k}^{1}(\mathbf{s})=F_{k-1}^{1}(\mathbf{f}_{k}(\mathbf{s}))f_{k}^{1}(\mathbf{s}).$
\end{lemma}
\proof For $k\ge 1,$ note that
\begin{equation}
Z_{k,p}=\sum_{j=1}^{Z_{k-1,1}+1}\xi_{p}^{j}(1,k)
+\sum_{j=1}^{Z_{k-1,2}}\xi_{p}^{j}(2,k)
+\sum_{j=1}^{Z_{k-1,3}}\xi_{p}^{j}(3,k),\ p\in \{1,2,3\}.\no
\end{equation}
Since for  $1\le p\le 3$ and $k\ge1,$ $\xi_{p}^{j}(m,k),j=1,2,\cdots$ are  independent copies of $\xi_{p}(m,k),$
we have
\begin{align}
F_{k}^{1}&(\mathbf{s})
=E_{\omega}\left[E_{\omega}\big(s_{1}^{Z_{k,1}}s_{2}^{Z_{k,2}}
s_{3}^{Z_{k,3}}\big|\mathbf{Z}_{k-1}\big)\right]\no\\
&= E_{\omega}\left(s_{1}^{\xi_{1}(1,k)}s_{2}^{\xi_{2}(1,k)}
s_{3}^{\xi_{3}(1,k)}\right)
E_{\omega}\bigg\{
\prod_{m=1}^{3}\big[E_{\omega}\big(s_{1}^{\xi_{1}(m,k)}
s_{2}^{\xi_{2}(m,k)}s_{3}^{\xi_{3}(m,k)}\big)\big]^{Z_{k-1,m}}\bigg\}\no\\
&=f_{k}^{1}(\mathbf{s})\cdot E_{\omega}\left\{[f_{k}^{1}(\mathbf{s})]^{Z_{k-1,1}}
[f_{k}^{2}(\mathbf{s})]^{Z_{k-1,2}}[f_{k}^{3}(\mathbf{s})]^{Z_{k-1,3}}\right\}\no\\
&=f_{k}^{1}(\mathbf{s})F_{k-1}^{1}(\mathbf{f}_{k}(\mathbf{s}))\no.
\end{align}
The lemma is proved. \qed

\begin{lemma}\label{Flem2}$F_{n}^{1}(\mathbf{s})
=\prod_{k=1}^{n}f_{k}^{1}(\mathbf{f}_{k+1}(\cdots\mathbf{f}_{n}(\mathbf{s})
\cdots)).$
\end{lemma}
\proof
With Lemma \ref{Flem1} in hands, we can prove Lemma \ref{Flem2}  by induction.\qed

\begin{lemma}\label{NEW} For $k\geq 1,$ we have
\begin{equation}\label{fks}
\mathbf{f}_{k}(\mathbf{s})
=\mathbf{1}-\frac{A_{k}(\mathbf{1}-\mathbf{s})}
{1+\mathbf{e}_{1}A_{k}(\mathbf{1}-\mathbf{s})}.\no
\end{equation}
\end{lemma}
\proof Using the formula $\frac{1}{(1-x)^{k+1}}=\sum_{m=0}^{+\infty}\frac{(m+k)!}{m!k!}x^{m}$, we have
\begin{align}
f_{k}^{1}(\mathbf{s})&=f_{k}^{3}(\mathbf{s})
=E_{\omega}(s_{1}^{\xi_{1}(1,k)}s_{2}^{\xi_{2}(1,k)}s_{3}^{\xi_{3}(1,k)})\no\\
&= \sum_{a=0}^{+\infty}\sum_{b=0}^{+\infty}\frac{(a+b)!}{a!b!}
(\omega_{k}(-1))^{a}(\omega_{k}(-2))^{b}\omega_{k}(1)s_{1}^{a}s_{2}^{b}\no\\
&=\frac{\omega_{k}(1)}{1-\omega_{k}(-1)s_{1}-\omega_{k}(-2)s_{2}}\no
\end{align}
and
\begin{align}
f_{k}^{2}(\mathbf{s})
&=E_{\omega}(s_{1}^{\xi_{1}(2,k)}s_{2}^{\xi_{2}(2,k)}s_{3}^{\xi_{3}(2,k)})\no\\
&= \sum_{a=0}^{+\infty}\sum_{b=0}^{+\infty}\frac{(a+b)!}{a!b!}
(\omega_{k}(-1))^{a}(\omega_{k}(-2))^{b}\omega_{k}(1)s_{1}^{a}s_{2}^{b}s_{3}\no\\
&=\frac{\omega_{k}(1)s_{3}}{1-\omega_{k}(-1)s_{1}-\omega_{k}(-2)s_{2}}.\no
\end{align}
Then some direct computation yields that
\begin{equation}
f_{k}^{1}(\mathbf{s})=f_{k}^{3}(\mathbf{s})
=1-\frac{\mathbf{\gamma}_{k}(\mathbf{1}-\mathbf{s})}
{1+\mathbf{\gamma}_{k}(\mathbf{1}-\mathbf{s})},
f_{k}^{2}(\mathbf{s})
=1-\frac{(a_{k},b_{k},1)(\mathbf{1}-\mathbf{s})}
{1+\mathbf{\gamma}_{k}(\mathbf{1}-\mathbf{s})},\no
\end{equation}
where $\mathbf{\gamma}_{k}=(a_{k},b_{k},0).$ The lemma follows. \qed

\begin{lemma}\label{Flem3} For $n\ge1$ and $1\le k\le n,$ we have
\begin{align}\label{mj1s}
\mathbf{f}_{1}(\mathbf{f}_{2}(\mathbf{f}_{3}(\cdots\mathbf{f}_{n}(\mathbf{s})
\cdots)))
&=\mathbf{1}-\frac{\prod_{j=1}^{n}A_{j}(\mathbf{1}-\mathbf{s})}
{1+\sum_{k=1}^{n}\mathbf{e}_{1}
\prod_{i=k}^{n}A_{i}(\mathbf{1}-\mathbf{s})},\no\\
f_{k}^{1}(\mathbf{f}_{k+1}(\cdots\mathbf{f}_{n}(\mathbf{s})\cdots))
&=1-\frac{\mathbf{e}_{1}\prod_{j=k}^{n}A_{j}(\mathbf{1}-\mathbf{s})}
{1+\sum_{j=k}^{n}\mathbf{e}_{1}
\prod_{i=j}^{n}A_{i}(\mathbf{1}-\mathbf{s})}.\no
\end{align}
\end{lemma}
\proof The lemma follows from Lemma \ref{NEW} by  induction.\qed

\noindent{\it Proof of Proposition \ref{Flem4}.} Taking Lemma \ref{Flem2} and Lemma \ref{Flem3}  together, we get
\begin{align*}
F_{n}^{1}(\mathbf{s})
=\prod_{k=1}^{n}\frac{1+\sum_{j=k+1}^{n}\mathbf{e}_{1}\prod_{i=j}^{n}A_{i}
(\mathbf{1}-\mathbf{s})}
{1+\sum_{j=k}^{n}\mathbf{e}_{1}\prod_{i=j}^{n}A_{i}(\mathbf{1}-\mathbf{s})}
=\frac{1}{1+\sum_{j=1}^{n}\mathbf{e}_{1}\prod_{i=j}^{n}A_{i}
(\mathbf{1}-\mathbf{s})},
\end{align*}
which implies (\ref{ff}). Moreover, since  $(\omega_1,\omega_2,...,\omega_n)$ and $(\omega_n,\omega_{n-1},...,\omega_1)$ share the same law under $\mathbb{P}^{\theta},$ it follows from (\ref{ff}) that $F_{n}^{1}(\mathbf{s})$ equals $\widetilde{F}_{n}^{1}(\mathbf{s})$ in $\mathbb P^\theta$-distribution. \qed

\begin{proposition}\label{fc}
Both $F_{n}^{1}(\mathbf{s})$ and $\widetilde{F}_{n}^{1}(\mathbf{s})$ converge in $\mathbb{P}^{\theta}$-distribution  to a p.g.f.
\begin{equation}
F^{1}(\mathbf{s})
:=\frac{\frac{1}{\hat{S}}}{1-\frac{\hat{S}^{1}}{\hat{S}}s_{1}
-\frac{\hat{S}^{2}}{\hat{S}}s_{2}
-\frac{\hat{S}^{3}}{\hat{S}}s_{3}},\no
\end{equation}
 where
\begin{equation}\label{si}
\hat{S}^{i}:=\mathbf{e}_{1}\big(A_{1}+A_{1}A_{2}+\cdots+A_{1}A_{2}\cdots A_{n}+\cdots\big)\mathbf{e}_{i}^{t},~i=1,2,3,
\end{equation} are all finite
and $\hat{S}:=1+\hat{S}^{1}+\hat{S}^{2}+\hat{S}^{3}$.
\end{proposition}
\proof Note that by (C2), $\gamma_A<0.$ Moreover, from (\ref{gab}), we have
\begin{equation}
\mathbb{P}^{\theta}\z(\exists N(\omega),\text{ s.t. }\forall n>N(\omega),
\|A_{0}A_{1}\cdots A_{n}\|\leq e^{n\gamma_A/2}\y)=1.\no
\end{equation}
Therefore $\mathbb{P}^{\theta}$-a.s., $\hat{S}^{i}$  defined in (\ref{si}) are all finite and
\begin{equation*}
\lim_{n\rightarrow\infty}\widetilde{S}_{n}^{i}
=\hat{S}^{i}, i=1,2,3 \text{ and } 
\lim_{n\rightarrow\infty}\widetilde{S}_{n}
=\hat{S}.
\end{equation*}
Consequently, $\widetilde{F}_{n}^{1}(\mathbf{s})$ converges $\mathbb P^\theta$-a.s. and hence in $\mathbb{P}^{\theta}$-distribution   to $F^{1}(\mathbf{s})$. Next, note that by Proposition \ref{Flem4},   $F_{n}^{1}(\mathbf{s})$ shares the same law with $\widetilde{F}_{n}^{1}(\mathbf{s}).$ Thus,  $F_{n}^{1}(\mathbf{s})$  converges also in $\mathbb{P}^{\theta}$-distribution to $F^{1}(\mathbf{s}).$
\qed
\subsection{$\{\mathbf Z_n\}$ as a Markov chain with a limit invariant distribution}
Since $\{\omega_{n}\}$ is an i.i.d. sequence under $\mathbb P^\theta,$ then  $\{\mathbf{Z}_{n}\}$ is a homogeneous Markov chain under $\mathbf P^{\theta},$ whose transition kernel can also be computed explicitly. This observation of Markov property ensures an application of the ergodic theorem for Markov chain to $\tilde\phi_\theta(\mathbf{Z}_n,\mathbf{Z}_{n+1}),~n\ge 0,$ see Proposition \ref{ergo} below.
\begin{proposition}\label{tran}
Under $\mathbf P^\theta,$ $\{\mathbf{Z}_{n}\}_{n\ge 0}$ is a homogeneous Markov chain, whose transition kernel $Q_{\theta}(\mathbf{x},\mathbf{y}),$  for any $\theta\in\Theta,$  is given by
\begin{align}\label{qk}
Q_{\theta}(\mathbf{x},\mathbf{y})
&=\frac{(x_{1}+x_{2}+x_{3}+y_{1}+y_{2})!}{y_{1}!y_{2}!(x_{1}+x_{2}+x_{3})!}\no\\
&\quad\quad\quad\quad\quad\times\iint\limits_{0\leq a_{1}+a_{2}\leq 1}a_{1}^{y_{1}}a_{2}^{y_{2}}
(1-a_{1}-a_{2})^{x_{1}+x_{2}+x_{3}+1}
d\nu_{\theta}(a_{1},a_{2}),
\end{align}
where $\mathbf{x}=(x_{1},x_{2},x_{3})\in\mathbb{N}^{3},$ $ \mathbf{y}=(y_{1},y_{2},y_{3})\in\mathbb{N}^{3}.$
\end{proposition}
\proof Let $\mathbf Z_{n}=(Z_{n,1},Z_{n,2},Z_{n,3})$ be individuals in the $n$-th generation conditioned on $\{\mathbf{Z}_{n-1}=(x_{1},x_{2},x_{3})\}$. Then for $p=1,2,3,$
\begin{equation}
Z_{n,p}=\sum_{j=1}^{x_{1}+1}\xi_{p}^{j}(1,n)
+\sum_{j=1}^{x_{2}}\xi_{p}^{j}(2,n)
+\sum_{j=1}^{x_{3}}\xi_{p}^{j}(3,n).\no
\end{equation}
The quenched p.g.f. of $\mathbf Z_{n},$ conditioned on $\{\mathbf{Z}_{n-1}=(x_{1},x_{2},x_{3})\},$ equals
\begin{eqnarray}
&&f_{\mathbf Z_{n}}(\mathbf{s})=[f_{n}^{1}(\mathbf{s})]^{x_{1}+1}
[f_{n}^{2}(\mathbf{s})]^{x_{2}}[f_{n}^{3}(\mathbf{s})]^{x_{3}}\no\\
&=&[1-\omega_{n}(-1)-\omega_{n}(-2)]^{x_{1}+x_{3}+1}
\frac{[(1-\omega_{n}(-1)-\omega_{n}(-2))s_{3}]^{x_{2}}}
{[1-\omega_{n}(-1)s_{1}-\omega_{n}(-2)s_{2}]^{x_{1}+x_{2}+x_{3}+1}}.\no
\end{eqnarray}
Some easy computation yields
\begin{align}
&\left.\frac{\partial ^{y_{1}+y_{2}+y_{3}}f_{\mathbf Z_{n}}(\mathbf{s})}{\partial y_{1}\partial y_{2}\partial y_{3}}\right|_{\mathbf{s}=(0,0,0)}=\frac{(x_{1}+x_{2}+x_{3}+y_{1}+y_{2})!y_{3}!}{(x_{1}+x_{2}+x_{3})!}\no\\
&\quad\quad\quad\quad\quad\quad\times
[\omega_{n}(-1)]^{y_{1}}[\omega_{n}(-2)]^{y_{2}}
[1-\omega_{n}(-1)-\omega_{n}(-2)]^{x_{1}+x_{2}+x_{3}+1}.\no
\end{align}
Therefore
\begin{align}
&P_{\omega}(\mathbf{Z}_{n}=(y_{1},y_{2},y_{3})|\mathbf{Z}_{n-1}=(x_{1},x_{2},x_{3}))
=\frac{(x_{1}+x_{2}+x_{3}+y_{1}+y_{2})!}{y_{1}!y_{2}!(x_{1}+x_{2}+x_{3})!}\no\\
&\quad\quad\quad\quad\quad\quad\times[\omega_{n}(-1)]^{y_{1}}[\omega_{n}(-2)]^{y_{2}}
[1-\omega_{n}(-1)-\omega_{n}(-2)]^{x_{1}+x_{2}+x_{3}+1},\no
\end{align}
and hence transition kernel $Q_{\theta}(\mathbf{x},\mathbf{y})$ coincides with the one defined in (\ref{qk}). \qed

\begin{theorem}\label{sd}
Assume  (C1) and (C2). Then for any $\theta\in\Theta,$  the chain $\{\mathbf{Z}_{n}\}_{n\ge0}$ is positive recurrent under $\mathbf P^\theta$ and admits a unique invariant probability measure $\pi_\theta$ satisfying
\begin{equation}
\lim_{n\rightarrow\infty}\mathbf{P}^{\theta}(\mathbf{Z}_{n}=(a,b,c))\label{ldv}
=\pi_{\theta}(a,b,c).
\end{equation}
 Moreover for all $(a,b,c)\in\mathbb{N}^{3},$
\begin{equation}
\pi_{\theta}(a,b,c)
=\mathbb{E}^{\theta}\bigg[\frac{(a+b+c)!}{a!b!c!}
\bigg(\frac{\hat{S}^{1}}{\hat{S}}\bigg)^{a}
\bigg(\frac{\hat{S}^{2}}{\hat{S}}\bigg)^{b}
\bigg(\frac{\hat{S}^{3}}{\hat{S}}\bigg)^{c}
\bigg(1-\frac{\hat{S}^{1}}{\hat{S}}-\frac{\hat{S}^{1}}{\hat{S}}
-\frac{\hat{S}^{3}}{\hat{S}}\bigg)\bigg].\no
\end{equation}
\end{theorem}
\begin{remark}
  As pointed out in \cite{cfl},  Key (\cite{key}, Theorem 3.3) showed that  $\mathbf{Z}_{n}$ converges in annealed law to a finite limit. Also, a construction of the limit was given in \cite{roi}. But here, since $\{\mathbf Z_n\}$  is Markovian and its offspring distributions are of the linear fractional type, we can compute explicitly the distribution of $\mathbf Z_n$ and its limit distribution as well.
\end{remark}

\proof Since for $i=1,2,3,$ $0<\frac{S_{n}^{i}}{S_{n}}<1,$ $\mathbb P^\theta$-a.s.,  thus by Proposition \ref{Flem4}, we see that, for any $(a,b,c)\in\mathbb{N}^{3},$
\begin{align*}
P_{\omega}(\mathbf{Z}_{n}=(a,b,c))
=\frac{(a+b+c)!}{a!b!c!}
\left(\frac{S_{n}^{1}}{S_{n}}\right)^{a}
\left(\frac{S_{n}^{2}}{S_{n}}\right)^{b}
\left(\frac{S_{n}^{3}}{S_{n}}\right)^{c}
\Big(1-\frac{S_{n}^{1}}{S_{n}}-\frac{S_{n}^{2}}{S_{n}}
-\frac{S_{n}^{3}}{S_{n}}\Big).
\end{align*}
As a result, using the fact
$(\omega_1,\omega_2,...,\omega_n)$ and $(\omega_n,\omega_{n-1},...,\omega_1)$ share the same law under $\mathbb{P}^{\theta}$ in the second equality,  we get
\begin{align}
\mathbf{P}^{\theta}&(\mathbf{Z}_{n}=(a,b,c))
=\mathbb{E}^{\theta}P_{\omega}(\mathbf{Z}_{n}=(a,b,c))\no\\
&=\mathbb{E}^{\theta}\bigg[\frac{(a+b+c)!}{a!b!c!}
\bigg(\frac{S_{n}^{1}}{S_{n}}\bigg)^{a}
\bigg(\frac{S_{n}^{2}}{S_{n}}\bigg)^{b}
\bigg(\frac{S_{n}^{3}}{S_{n}}\bigg)^{c}
\bigg(1-\frac{S_{n}^{1}}{S_{n}}-\frac{S_{n}^{2}}{S_{n}}
-\frac{S_{n}^{3}}{S_{n}}\bigg)\bigg]\no\\
&=\mathbb{E}^{\theta}\bigg[\frac{(a+b+c)!}{a!b!c!}
\bigg(\frac{\widetilde{S}_{n}^{1}}{\widetilde{S}_{n}}\bigg)^{a}
\bigg(\frac{\widetilde{S}_{n}^{2}}{\widetilde{S}_{n}}\bigg)^{b}
\bigg(\frac{\widetilde{S}_{n}^{3}}{\widetilde{S}_{n}}\bigg)^{c}
\bigg(1-\frac{\widetilde{S}_{n}^{1}}{\widetilde{S}_{n}}
-\frac{\widetilde{S}_{n}^{2}}{\widetilde{S}_{n}}
-\frac{\widetilde{S}_{n}^{3}}{\widetilde{S}_{n}}\bigg)\bigg].\no
\end{align}
On the other hand, since $0<\frac{\widetilde{S}_{n}^{i}}{\widetilde{S}_{n}}
<1,$ ~$i=1,2,3,$  and $0<\frac{\widetilde{S}_{n}^{1}}{\widetilde{S}_{n}}
+\frac{\widetilde{S}_{n}^{2}}{\widetilde{S}_{n}}
+\frac{\widetilde{S}_{n}^{3}}{\widetilde{S}_{n}}<1,$  $\mathbb{P}^{\theta}$-a.s., then by the dominated convergence theorem, we have
\begin{align}
\lim_{n\rightarrow\infty}&\mathbf{P}^{\theta}(\mathbf{Z}_{n}=(a,b,c))\no\\
&=\mathbb{E}^{\theta}\bigg[\frac{(a+b+c)!}{a!b!c!}
\bigg(\frac{\hat{S}^{1}}{\hat{S}}\bigg)^{a}
\bigg(\frac{\hat{S}^{2}}{\hat{S}}\bigg)^{b}
\bigg(\frac{\hat{S}^{3}}{\hat{S}}\bigg)^{c}
\bigg(1-\frac{\hat{S}^{1}}{\hat{S}}-\frac{\hat{S}^{1}}{\hat{S}}
-\frac{\hat{S}^{3}}{\hat{S}}\bigg)\bigg]\no\\
&=:\pi_{\theta}(a,b,c).\no
\end{align}
Moreover, by Proposition \ref{fc}, we have $0<\frac{\hat{S}^{i}}{\hat{S}}<1,$ $i=1,2,3,$ $\mathbb{P}^{\theta}$-a.s. Thus, by Fubini-Tonelli's theorem,  we get $\sum_{a,b,c\geq0}\pi_{\theta}(a,b,c)=1$. So  $\pi_{\theta}$ is a probability measure on $\mathbb{N}^{3}.$  Due to (\ref{ldv}), $\pi_\theta$ is invariant.  Moreover, since the  kernel $Q_{\theta}(\mathbf{x},\mathbf{y})$ is positive and  $\nu_{\theta}$ is not degenerate, the chain $\{\mathbf{Z}_{n}\}$ is irreducible. Therefore, $\{\mathbf Z_n\}$ is positive recurrent and $\pi_{\theta}$ is unique. This completes the proof.
\qed

\section{Consistency}\label{consi}

In this section, we study the consistency of the estimator. Under certain conditions, we  show that $l_{n}(\theta)/n$ converges  to a finite limit $l(\theta),$ which identifies the true parameter value. Then, the consistency of the estimator can be established by some standard statistical arguments.

To begin with, let $\{\mathbf{\widetilde{Z}}_{n}\}_{n\ge0}$ be a process
with transition kernel $Q_{\theta^{\ast}}(\mathbf{x},\mathbf{y})$ and initial distribution $\pi_{\theta^{\ast}}.$ Then
$\{\mathbf{\widetilde{Z}}_{n}\}$ is a
stationary Markov chain under $\mathbf P^*.$
 We have shown that under $\mathbf{P}^{\ast},$ $\{\mathbf{Z}_{n}\}$ is an irreducible Markov chain   with transition kernel $Q_{\theta^{\ast}}(\mathbf{x},\mathbf{y}).$   We have the following ergodic theorem.

\begin{proposition}\label{ergo} Assume (C1) and (C2) hold. Then
for all $\theta\in\Theta,$  we have
\begin{equation}\label{ergo1}\no
\lim_{n\rightarrow\infty}\frac{1}{n}\sum_{k=0}^{n-1}
\tilde{\phi}_{\theta}(\mathbf{Z}_{k},
\mathbf{Z}_{k+1})
=\mathbf{E}^{\ast}[\tilde{\phi}_{\theta}
(\mathbf{\widetilde{Z}}_{0},\mathbf{\widetilde{Z}}_{1})],
~\mathbf{P}^{\ast}\text{-a.s.}
\end{equation}
\end{proposition}

\proof Note that  $-\infty<\tilde{\phi}_{\theta}
(\mathbf{\widetilde{Z}}_{0},\mathbf{\widetilde{Z}}_{1})\le 0.$
If $\mathbf{E}^{\ast}[\tilde{\phi}_{\theta}
(\mathbf{\widetilde{Z}}_{0},\mathbf{\widetilde{Z}}_{1})]>-\infty,$ then the proposition is the same as Proposition 2.4 in \cite{FL14}.
If  $\mathbf{E}^{\ast}[\tilde{\phi}_{\theta}
(\mathbf{\widetilde{Z}}_{0},\mathbf{\widetilde{Z}}_{1})]=-\infty,$ then for $M<0,$ we have $\mathbf P^*$-a.s.,
\begin{align*}
\limsup_{n\rto}\frac{1}{n}\sum_{k=0}^{n-1}
\tilde{\phi}_{\theta}(\mathbf{Z}_{k},\mathbf{Z}_{k+1})
&\le \limsup_{n\rightarrow\infty}
\frac{1}{n}\sum_{k=0}^{n-1}
(\tilde{\phi}_{\theta}(\mathbf{Z}_{k},\mathbf{Z}_{k+1})\vee M)\no\\
&=\mathbf{E}^{\ast}[\tilde{\phi}_{\theta}
(\mathbf{\widetilde{Z}}_{0},\mathbf{\widetilde{Z}}_{1})\vee M].\no
\end{align*}
Letting $M\rightarrow-\infty$ we get $\mathbf P^*$-a.s.,
$$\limsup_{n\rto}\frac{1}{n}\sum_{k=0}^{n-1}
\tilde{\phi}_{\theta}(\mathbf{Z}_{k},\mathbf{Z}_{k+1})\le \mathbf{E}^{\ast}[\tilde{\phi}_{\theta}
(\mathbf{\widetilde{Z}}_{0},\mathbf{\widetilde{Z}}_{1})]=-\infty.$$
The proposition is proved. \qed

Next let
\begin{equation}
l(\theta):=\mathbf{E}^{\ast}[\tilde{\phi}_{\theta}
(\mathbf{\widetilde{Z}}_{0},\mathbf{\widetilde{Z}}_{1})].\no
\end{equation}
\begin{lemma}\label{finite}
Assume (C1)-(C3) hold. Then   $l(\theta)$ is finite for any  $\theta\in\Theta.$
\end{lemma}
\proof By Jensen's inequality, for all $\mathbf x=(x_{1},x_{2},x_{3})\in \mathbb{N}^{3},$~$(y_{1},y_{2},y_{3})\in \mathbb{N}^{3},$
\begin{align}
\log&\iint_{0\leq a_{1}+a_{2}\leq 1}a_{1}^{y_{1}}a_{2}^{y_{2}}
(1-a_{1}-a_{2})^{x_{1}+x_{2}+x_{3}+1}
d\nu_{\theta}(a_{1},a_{2})\no\\
&\geq
y_{1}\mathbb{E}^{\theta}[\log(\omega_{0}(-1))]
+y_{2}\mathbb{E}^{\theta}[\log(\omega_{0}(-2))]
+(|\mathbf x|+1)\mathbb{E}^{\theta}[\log(\omega_{0}(1))].\no
\end{align}
Thus
\begin{align}\label{lower}
\frac{1}{n}\sum_{k=0}^{n-1}\tilde{\phi}_{\theta}
(\mathbf{Z}_{k},\mathbf{Z}_{k+1})
&\geq
\frac{1}{n}\sum_{k=0}^{n-1}\mathbf{Z}_{k+1}
\cdot(\mathbb{E}^{\theta}[\log(\omega_{0}(-1))],
\mathbb{E}^{\theta}[\log(\omega_{0}(-2))],0)^{t}\no\\
&\quad\quad+\frac{1}{n}\sum_{k=0}^{n-1}(|\mathbf{Z}_{k}|+1)
\cdot\mathbb{E}^{\theta}[\log(\omega_{0}(1))].
\end{align}
In view of Corollary \ref{mm}, (\ref{uz})  and the stationarity of the environment, we get
\begin{align}\label{ze}
\mathbf{E}^{\ast}(\tilde{\mathbf{Z}}_{0})
&=\lim_{n\rightarrow\infty}\mathbf{E}^{\ast}(\mathbf{Z}_{n})
=\lim_{n\rto}\mathbb{E}^{\ast}\Big(\sum_{k=1}^{n}
\mathbf{e}_{1}A_{k}A_{k+1}\cdots A_{n}\Big)\no\\
&=\lim_{n\rto}\mathbb{E}^{\ast}\Big(\sum_{k=1}^{n}
\mathbf{e}_{1}A_{k}A_{k-1}\cdots A_{1}\Big).
\end{align}
By the ergodic theorem of Markov chain(Norris \cite{mc}, Theorem 1.10.2),  it follows that $\mathbf{P}^{\ast}$-a.s.,
\begin{equation}
\lim_{n\rightarrow+\infty}\frac{1}{n}\sum_{k=0}^{n-1}\mathbf{Z}_{k+1}
=\mathbf{E}^{\ast}(\tilde{\mathbf{Z}}_{0})\text{ and }
\lim_{n\rightarrow+\infty}\frac{1}{n}\sum_{k=0}^{n-1}(|\mathbf{Z}_{k}|+1)
=\mathbf{E}^{\ast}(|\tilde{\mathbf{Z}}_{0}|)+1.\no
\end{equation}
Consequently,  by (\ref{ze}), (C1) and (C3), the limit of the lower bound in (\ref{lower}) is finite and so is $l(\theta)$.\qed

%
%

Note that by (\ref{lz}), we have
$
l_{n}(\theta)=
\sum_{k=0}^{n-1}\tilde{\phi}_{\theta}
(\mathbf{Z}_{k},\mathbf{Z}_{k+1}),
$  in $\mathbf{P}^{\ast}$-distribution.
Then it follows from Proposition \ref{ergo} that
\begin{equation}
\lim_{n\rightarrow\infty}\frac{l_{n}(\theta)}{n}
=l(\theta),~\text{ in }\mathbf{P}^{\ast}\text{-probability}.\no
\end{equation}
To sum up, we have the following theorem.
\begin{theorem}\label{thm4.1} Suppose that (C1)-(C3) hold. Then, there exists a finite deterministic limit $l(\theta)$ such that
$
\lim_{n\rightarrow+\infty}\frac{l_{n}(\theta)}{n}=l(\theta),
$ in $\mathbf{P}^{\ast}$-probability.
\end{theorem}
\begin{remark}
Indeed, for the convergence in Theorem \ref{thm4.1}, we can say more. Suppose (C1)-(C3) hold. Fix an open $U\subset \Theta$ and  define $\widetilde{\Phi}_{U}=\sup_{\theta\in U}\tilde{\phi}_{\theta}$. Then
$\frac{1}{n}\sum_{x=0}^{n-1}\sup_{\theta\in U}
\phi_{\theta}(\mathbf{L}_{x}^{n},\mathbf{L}_{x+1}^{n},\mathbf{L}_{x+2}^{n})=\frac{1}{n}\sum_{k=0}^{n-1}\widetilde{\Phi}_{U}
(\mathbf{Z}_{k},\mathbf{Z}_{k+1})$ in $\mathbf{P}^{\ast}$-distribution.
Since $\widetilde{\Phi}_{U}$ is always negative, $\mathbf{E}^{\ast}[\widetilde{\Phi}_{U}(\tilde{\mathbf{Z}}_{0},\tilde{\mathbf{Z}}_{1})]$ is well defined. Thus by the ergodic theorem in Proposition \ref{ergo}, we conclude that $$
\lim_{n\rightarrow+\infty}\frac{1}{n}\sum_{x=0}^{n-1}
\displaystyle\sup_{\theta\in U}
\phi_{\theta}(\mathbf{L}_{x}^{n},\mathbf{L}_{x+1}^{n},\mathbf{L}_{x+2}^{n})
=\mathbf{E}^{\ast}\Big[\sup_{\theta\in U}\tilde{\phi}_{\theta}
(\mathbf{\widetilde{Z}}_{0},\mathbf{\widetilde{Z}}_{1})\Big],
$$ in $\mathbf{P}^{\ast}$-probability. That is, $l_n(\theta)/n$ converges to $l(\theta)$ in $\mathbf{P}^{\ast}$-probability in a `locally unform' manner as $n\rto. $\end{remark}

Similar to \cite{cfl}, we have from the following theorem  that  $l(\theta)$ identifies the true value $\theta^{\ast}$  as the unique point where it attains its maximum.
\begin{theorem}\label{prop4.2}
Suppose that (C1)-(C5) hold. Then for any fixed $\varepsilon>0,$
\begin{equation}
\displaystyle\sup_{\theta:~||\theta-\theta^{\ast}||\geq\varepsilon} l(\theta)<l(\theta^{\ast}).\no
\end{equation}
\end{theorem}
{\noindent\it Proof of Theorem \ref{consist}.} Since  $\Theta$ is compact, by applying some classical arguments of  consistency for $M$-estimator(see \cite{van}, Section 5.2.1 therein), Theorem \ref{consist} follows from Theorem \ref{thm4.1} and Theorem \ref{prop4.2}. \qed

\end{document}